\documentclass[reqno]{amsart}%
\usepackage{amssymb}
\usepackage{amsfonts}%
\usepackage{amsmath}%
\usepackage{graphicx}

\providecommand{\U}[1]{\protect \rule{.1in}{.1in}}
\newtheorem{theorem}{Theorem}
\theoremstyle{plain}

\numberwithin{equation}{section}

\begin{document}
	\title[Derivatives of Humbert confluent hypergeometric functions]{ Derivatives of Humbert confluent hypergeometric functions with respect to their parameters}
	\author[A. Shehata, R. \c{S}ahin ,  O. Ya\u{g}c\i, S.I. Moustafa ]{Ayman Shehata$^{1}$, Recep \c{S}ahin$^{2}$, O\u{g}uz Ya\u{g}c\i$^{2}$, Shimaa I. Moustafa$^{3}$}
	
	
\begin{abstract}
	Humbert confluent hypergeometric functions of two variables arise in many problems of mathematical physics and applied analysis, yet their behaviour with respect to parameters has not been systematically studied. In this paper we investigate derivatives with respect to numerator and denominator parameters for the seven classical Humbert functions
	\(\Phi_{1}\), \(\Phi_{2}\), \(\Phi_{3}\), \(\Psi_{1}\), \(\Psi_{2}\), \(\Xi_{1}\) and \(\Xi_{2}\).
	Using their double–series representations together with elementary properties of the Gamma and digamma functions, we derive explicit formulas for first–order parameter derivatives and express them in compact form in terms of Srivastava’s triple hypergeometric function \(F^{(3)}\). By differentiating the underlying partial differential equations, we further obtain simple operator recurrences for derivatives of arbitrary order, which yield closed differentiation and reduction formulas in terms of contiguous Humbert functions. Finally, we indicate how these results lead to Taylor-type parameter expansions and illustrate their use with basic numerical examples and plots.
\end{abstract}

	\maketitle
	
\section{Introduction}

Hypergeometric functions and their multivariable analogues play a central rôle
in the theory of special functions and in many areas of mathematical physics,
engineering and applied analysis. Starting from the classical Gauss and Kummer
functions of one variable, various generalizations have been introduced,
including the Appell and Lauricella families of two or more variables (see,
for example, \cite{appell1926fonctions,srivastava1985multiple}). Among their
confluent limits, the seven functions introduced by P.~Humbert form a
distinguished class of two–variable confluent hypergeometric functions, now
customarily denoted by
\(\Phi_{1}\), \(\Phi_{2}\), \(\Phi_{3}\),
\(\Psi_{1}\), \(\Psi_{2}\),
\(\Xi_{1}\) and \(\Xi_{2}\).
They admit simple double–series representations in terms of Pochhammer symbols
and Gamma functions and satisfy systems of linear partial differential
equations with polynomial coefficients.

In many applications, the parameters of these functions carry direct physical,
geometric or probabilistic meaning, and one is interested not only in the
functions themselves, but also in their variation with respect to the
parameters. Derivatives with respect to numerator or denominator parameters
arise naturally in sensitivity analysis, perturbation methods, analytic
continuation and the derivation of asymptotic expansions. For single–variable
hypergeometric functions such as \({}_{2}F_{1}\), \({}_{1}F_{1}\) and more
general \({}_{p}F_{q}\) series, parameter derivatives have been studied
extensively and can be expressed in terms of polygamma functions and shifted
hypergeometric functions (see, e.g., \cite{exton1976multiple} and the
references therein). By contrast, a systematic treatment of parameter
derivatives for the Humbert confluent hypergeometric functions of two variables
has been lacking.

The aim of this work is to provide such a systematic treatment. Starting from
the double–series definitions of the seven Humbert functions, and using
elementary properties of the Gamma, digamma and polygamma functions together
with standard identities for the Pochhammer symbol, we first derive explicit
formulas for the derivatives with respect to each numerator and denominator
parameter. These first–order parameter derivatives are then recast in a compact
and unified way in terms of Srivastava’s triple hypergeometric function
\(F^{(3)}\), which plays a natural rôle as a building block for multivariable
parameter differentiation.

A second ingredient of our approach is the use of the systems of partial
differential equations satisfied by the Humbert functions. By differentiating
these PDEs with respect to the parameters, we obtain simple operator recurrence
relations which express \(n\)th–order parameter derivatives in terms of
lower–order ones. Combining these recurrences with shift identities for
Pochhammer symbols leads to closed differentiation and reduction formulas
relating parameter derivatives to contiguous Humbert functions and to higher–
order derivatives with respect to the variables.

Finally, to illustrate the applicability of the theoretical results, we present
a short numerical study for selected Humbert functions. In particular, we
compute sample values and two– and three–dimensional plots for \(\Phi_{1}\) and
its derivative with respect to a numerator parameter, using the series
representations implied by our formulas. These numerical illustrations confirm
that the parameter–derivative calculus developed in this paper can be
implemented efficiently and provides a practical tool for applications.

The paper is organized as follows. In Section~\ref{sec:preliminaries} we recall
the basic notation and properties of the Gamma and polygamma functions, the
Pochhammer symbol, the Humbert confluent hypergeometric functions and
Srivastava’s triple hypergeometric function \(F^{(3)}\). In
Section~\ref{sec:nth-order} we derive general recurrence relations for
\(n\)th–order derivatives with respect to the parameters by differentiating
the underlying systems of partial differential equations. Section
\ref{sec:applications} contains explicit differentiation and reduction formulas
in terms of contiguous Humbert functions. Numerical examples and graphical
illustrations are presented in Section~\ref{sec:numerical}. Concluding remarks
and some perspectives for further work are given in
Section~\ref{sec:conclusion}.

	\section{Preliminaries}
	\label{sec:preliminaries}
	
	In this section we collect the basic notation and auxiliary results which will be used throughout the paper. Unless otherwise stated, all parameters are complex and chosen so as to avoid poles of the Gamma function, and the variables $x$ and $y$ are complex numbers lying in the domains of convergence explicitly indicated below.
	
	\subsection{Gamma function, Pochhammer symbol and polygamma functions}
	
	We recall that the Euler Gamma function $\Gamma(z)$ is defined for $\Re(z)>0$ by
	\[
	\Gamma(z) = \int_{0}^{\infty} t^{z-1} e^{-t}\, dt,
	\]
	and is extended to a meromorphic function on $\mathbb{C}$ with simple poles at the non–positive integers. The rising factorial or Pochhammer symbol $(a)_n$ is given by
	\[
	(a)_0 := 1, \qquad
	(a)_n := a(a+1)\cdots(a+n-1) = \frac{\Gamma(a+n)}{\Gamma(a)}, \qquad n \in \mathbb{N}.
	\]
	
	The logarithmic derivative of the Gamma function is the digamma function
	\begin{equation}\label{eq:digamma-def}
		\Psi(z) := \frac{d}{dz} \log \Gamma(z)
		= \frac{\Gamma'(z)}{\Gamma(z)} ,
	\end{equation}
	and its higher derivatives
	\[
	\Psi_r(z) := \frac{d^r}{dz^r}\Psi(z), \qquad r \in \mathbb{N},
	\]
	are called the polygamma functions. Using the well–known representation of $\Psi$ in terms of a convergent series, one obtains the identity
	\begin{equation}\label{eq:digamma-shift}
		\Psi(z+n) - \Psi(z)
		= \sum_{k=0}^{n-1} \frac{1}{z+k},
		\qquad n \in \mathbb{N},
	\end{equation}
	and, more generally,
	\begin{equation}\label{eq:polygamma-shift}
		\Psi_r(z+n) - \Psi_r(z)
		= (-1)^r r! \sum_{k=0}^{n-1} \frac{1}{(z+k)^{r+1}},
		\qquad r \in \mathbb{N}_0,\; n \in \mathbb{N}.
	\end{equation}
	
	From \eqref{eq:digamma-def}–\eqref{eq:digamma-shift} and the representation of $(z)_n$ in terms of Gamma functions, we readily obtain the derivative of the Pochhammer symbol with respect to its parameter:
	\begin{equation}\label{eq:pochhammer-derivative}
		\frac{d}{dz}(z)_n
		= (z)_n\bigl[\Psi(z+n)-\Psi(z)\bigr]
		= (z)_n \sum_{k=0}^{n-1} \frac{1}{z+k},
		\qquad n \in \mathbb{N}.
	\end{equation}
	Similarly, differentiation of the reciprocal of a Pochhammer symbol yields
	\begin{equation}\label{eq:denominator-derivative}
		\frac{d}{dz} \frac{1}{(z)_n}
		= - \frac{1}{(z)_n} \sum_{k=0}^{n-1} \frac{1}{z+k},
		\qquad n \in \mathbb{N},
	\end{equation}
	which will be used below for derivatives with respect to denominator parameters.
	
	Later on we shall also employ a simple rearrangement formula for double series of the form
	\begin{equation}\label{eq:double-sum-rearrange}
		\sum_{n=0}^{\infty} \sum_{k=0}^{\infty} A(k,n)
		= \sum_{n=0}^{\infty} \sum_{k=0}^{n} A\bigl(k, n-k\bigr),
	\end{equation}
	valid whenever both sides converge absolutely. This identity allows us to convert sums over independent indices into sums over triangular regions, which will be convenient when expressing derivatives in terms of triple hypergeometric series.
	
	\subsection{Humbert confluent hypergeometric functions of two variables}
	
	We now recall the seven Humbert confluent hypergeometric functions of two variables, which are confluent forms of the classical Appell functions. They are defined in terms of double power series as follows:
	\begin{align}
		\Phi_{1}(a,b;c;x,y)
		&:= \sum_{m,n=0}^{\infty}
		\frac{(a)_{m+n} (b)_m}{(c)_{m+n}\, m!\, n!}\,
		x^m y^n,
		& |x| < 1,\; & |y| < \infty, \label{eq:Phi1-def} \\[0.4em]
		\Phi_{2}(a,b;c;x,y)
		&:= \sum_{m,n=0}^{\infty}
		\frac{(a)_m (b)_n}{(c)_{m+n}\, m!\, n!}\,
		x^m y^n,
		& |x| < \infty,\; & |y| < \infty, \label{eq:Phi2-def} \\[0.4em]
		\Phi_{3}(a;b;x,y)
		&:= \sum_{m,n=0}^{\infty}
		\frac{(a)_m}{(b)_{m+n}\, m!\, n!}\,
		x^m y^n,
		& |x| < \infty,\; & |y| < \infty. \label{eq:Phi3-def}
	\end{align}
	The two Humbert functions of $\Psi$–type are given by
	\begin{align}
		\Psi_{1}(a,b;c,d;x,y)
		&:= \sum_{m,n=0}^{\infty}
		\frac{(a)_{m+n} (b)_m}{(c)_m (d)_n\, m!\, n!}\,
		x^m y^n,
		& |x| < 1,\; & |y| < \infty, \label{eq:Psi1-def} \\[0.4em]
		\Psi_{2}(a;b,c;x,y)
		&:= \sum_{m,n=0}^{\infty}
		\frac{(a)_{m+n}}{(b)_m (c)_n\, m!\, n!}\,
		x^m y^n,
		& |x| < \infty,\; & |y| < \infty, \label{eq:Psi2-def}
	\end{align}
	and the two Humbert functions of $\Xi$–type are defined by
	\begin{align}
		\Xi_{1}(a,b,c;d;x,y)
		&:= \sum_{m,n=0}^{\infty}
		\frac{(a)_m (b)_n (c)_m}{(d)_{m+n}\, m!\, n!}\,
		x^m y^n,
		& |x| < 1,\; & |y| < \infty, \label{eq:Xi1-def} \\[0.4em]
		\Xi_{2}(a,b;c;x,y)
		&:= \sum_{m,n=0}^{\infty}
		\frac{(a)_m (b)_m}{(c)_{m+n}\, m!\, n!}\,
		x^m y^n,
		& |x| < 1,\; & |y| < \infty. \label{eq:Xi2-def}
	\end{align}
	
	Each of these double series defines an analytic function in the indicated domain of convergence, and admits analytic continuation in the parameters $a,b,c,d$ provided that poles of the Gamma function are avoided. The functions \eqref{eq:Phi1-def}–\eqref{eq:Xi2-def} satisfy systems of linear partial differential equations with polynomial coefficients; these systems will later be used to derive recurrence relations for higher–order derivatives with respect to the parameters.
	
	\subsection{Srivastava's triple hypergeometric function}
	
	A central rôle in our analysis is played by Srivastava's triple hypergeometric function, which provides a natural receptacle for the expressions arising from parameter differentiation of the Humbert functions. Following the notation used in the literature, we write
	\begin{equation}\label{eq:F3-def}
		F^{(3)}\!\left[
		\begin{array}{c}
			(a) :: (b); (b'); (b'') : (c); (c'); (c'') \\
			(e) :: (g); (g'); (g'') : (h); (h'); (h'')
		\end{array}
		\,\bigg|\, x,y,z
		\right]
		:= \sum_{m,n,p=0}^{\infty}
		\frac{A(m,n,p)}{m!\,n!\,p!}\,
		x^m y^n z^p,
	\end{equation}
	where the coefficient $A(m,n,p)$ is given by
	\begin{equation}\label{eq:F3-coeff}
		A(m,n,p)
		= \frac{
			\displaystyle
			\prod_{i=1}^{A} (a_i)_{m+n+p}
			\prod_{i=1}^{B} (b_i)_{m+n}
			\prod_{i=1}^{B'} (b'_i)_{n+p}
			\prod_{i=1}^{B''} (b''_i)_{m+p}
			\prod_{i=1}^{C} (c_i)_{m}
			\prod_{i=1}^{C'} (c'_i)_{n}
			\prod_{i=1}^{C''} (c''_i)_{p}
		}{
			\displaystyle
			\prod_{i=1}^{E} (e_i)_{m+n+p}
			\prod_{i=1}^{G} (g_i)_{m+n}
			\prod_{i=1}^{G'} (g'_i)_{n+p}
			\prod_{i=1}^{G''} (g''_i)_{m+p}
			\prod_{i=1}^{H} (h_i)_{m}
			\prod_{i=1}^{H'} (h'_i)_{n}
			\prod_{i=1}^{H''} (h''_i)_{p}
		}.
	\end{equation}
Here $(a)$ denotes the collection of parameters $a_1,\dots,a_A$, and similarly for the other grouped parameters $(b)$, $(b')$, $(b'')$, $(c)$, $(c')$, $(c'')$, $(e)$, $(g)$, $(g')$, $(g'')$, $(h)$, $(h')$ and $(h'')$. The precise conditions for convergence of the triple series \eqref{eq:F3-def} can be found in the standard references on multiple hypergeometric functions and will not be repeated here; in all subsequent applications, $(x,y,z)$ will be chosen so that the corresponding series is absolutely convergent.
	
In the sequel we shall make use of several specializations of \eqref{eq:F3-def} in which many of the parameter groups are empty or contain only a single parameter. In particular, the first–order parameter derivatives of the Humbert functions in Section~\ref{sec:nth-order} will be expressed in terms of $F^{(3)}$ with carefully chosen parameter arrays, while higher–order derivatives will be related to iterated differential operators acting on such triple hypergeometric series.
	
\section{Nth-order derivatives of Humbert confluent hypergeometric functions
	with respect to parameters}
\label{sec:nth-order}

In this section we derive recursive formulas for derivatives of arbitrary order
with respect to the parameters of the Humbert confluent hypergeometric
functions. The key idea is to make systematic use of the linear partial
differential equations (PDEs) satisfied by each Humbert function and to
differentiate these PDEs with respect to the parameters. Since the operators
involved are linear in the parameters, this procedure leads to simple
recurrence relations which express the $n$th parameter–derivative in terms of
the $(n-1)$st one.

Throughout this section we use the shorthand notation
\[
p := \frac{\partial}{\partial x}, \qquad
q := \frac{\partial}{\partial y}, \qquad
r := \frac{\partial^2}{\partial x^2}, \qquad
s := \frac{\partial^2}{\partial x\,\partial y}, \qquad
t := \frac{\partial^2}{\partial y^2}.
\]

\subsection{The case of $\Phi_{1}$}

We begin with the function
\[
\Phi_{1}(a,b;c;x,y)
= \sum_{m,n=0}^{\infty}
\frac{(a)_{m+n}(b)_m}{(c)_{m+n}\,m!\,n!}\,x^m y^n.
\]

It is known that $\Phi_{1}$ satisfies a pair of second–order linear PDEs of
hypergeometric type; see, for example, \cite{appell1926fonctions,srivastava1985multiple}.
For our purposes it is convenient to write these equations in the compact
operator form
\begin{align}
	D \Phi_{1}(a,b;c;x,y) &= 0, \label{eq:Dphi1}\\
	M \Phi_{1}(a,b;c;x,y) &= 0, \label{eq:Mphi1}
\end{align}
where
\begin{align}
	D
	&= x(1-x)\,r + y(1-x)\,s + y(1-y)\,t
	+ \bigl[c - (a+b+1)x\bigr]\,p
	- b y\,q - ab,
	\label{eq:D-operator}\\[0.4em]
	M
	&= y\,t + x\,s + (c-y)\,q - x\,p - a.
	\label{eq:M-operator}
\end{align}
Here $D$ and $M$ act on functions of $(x,y)$, but their coefficients depend
linearly on the parameters $a,b,c$. This dependence is what allows us to
generate parameter–derivative relations in a straightforward way.

\subsubsection*{First–order parameter derivatives}

To illustrate the method, consider first the derivative with respect to $a$.
Differentiating \eqref{eq:Dphi1} with respect to $a$ and using the fact that
$\partial/\partial a$ commutes with $D$ except through the coefficients, we
obtain
\[
D\!\left(\frac{\partial \Phi_{1}}{\partial a}\right)
+ \frac{\partial D}{\partial a}\,\Phi_{1} = 0.
\]
From \eqref{eq:D-operator} we compute
\[
\frac{\partial D}{\partial a}
= -x\,p - b,
\]
so that
\begin{equation}\label{eq:first-a}
	D\!\left(\frac{\partial \Phi_{1}}{\partial a}\right)
	= (x\,p + b)\,\Phi_{1}.
\end{equation}
The right–hand side is particularly simple: it consists of a first–order
differential operator in $(x,y)$ applied to the original function $\Phi_{1}$.
Using similar calculations for the derivatives with respect to $b$ and $c$ we
obtain
\begin{align}
	D\!\left(\frac{\partial \Phi_{1}}{\partial b}\right)
	&= (x\,p + y\,q + a)\,\Phi_{1}, \label{eq:first-b}\\
	D\!\left(\frac{\partial \Phi_{1}}{\partial c}\right)
	&= -p\,\Phi_{1}. \label{eq:first-c}
\end{align}
On the other hand, differentiating the second PDE \eqref{eq:Mphi1} with respect
to $a,b,c$ yields
\begin{align}
	M\!\left(\frac{\partial \Phi_{1}}{\partial a}\right)
	&= \Phi_{1}, \label{eq:first-a-M}\\
	M\!\left(\frac{\partial \Phi_{1}}{\partial b}\right)
	&= 0, \label{eq:first-b-M}\\
	M\!\left(\frac{\partial \Phi_{1}}{\partial c}\right)
	&= -q\,\Phi_{1}. \label{eq:first-c-M}
\end{align}
Equations \eqref{eq:first-a}–\eqref{eq:first-c-M} are the basic relations for
the first–order parameter derivatives of $\Phi_{1}$ obtained purely from the
PDEs.

In many applications it is convenient to rewrite the right–hand sides in terms
of contiguous Humbert functions such as $\Phi_{1}(a+1,b;c;x,y)$ and
$\Phi_{1}(a,b+1;c;x,y)$. Such representations will be used later in connection
with the differentiation formulas of Section~\ref{sec:applications}; for the
moment we keep the simpler operator form, which is sufficient to obtain
recurrence relations for higher–order derivatives. 

\subsubsection*{Recursive formulas for the $n$th derivatives}

We now differentiate the relations
\eqref{eq:first-a}–\eqref{eq:first-c} repeatedly with respect to the parameters.
For instance, applying $\partial^{n-1}/\partial a^{n-1}$ to
\eqref{eq:first-a} and using linearity of $D$ gives
\[
D\!\left(\frac{\partial^{n} \Phi_{1}}{\partial a^{n}}\right)
= (x\,p + b)\,\frac{\partial^{n-1}\Phi_{1}}{\partial a^{n-1}}
+ (n-1)\,\frac{\partial}{\partial a}(x\,p + b)\,
\frac{\partial^{n-2}\Phi_{1}}{\partial a^{n-2}}.
\]
Since the coefficient $x\,p + b$ is independent of $a$, its derivative with
respect to $a$ vanishes and the second term disappears. Thus we simply have
\[
D\!\left(\frac{\partial^{n} \Phi_{1}}{\partial a^{n}}\right)
= (x\,p + b)\,\frac{\partial^{n-1}\Phi_{1}}{\partial a^{n-1}}.
\]
By the same reasoning, repeated differentiation of \eqref{eq:first-b} and
\eqref{eq:first-c} yields
\[
D\!\left(\frac{\partial^{n} \Phi_{1}}{\partial b^{n}}\right)
= (x\,p + y\,q + a)\,\frac{\partial^{n-1}\Phi_{1}}{\partial b^{n-1}},
\qquad
D\!\left(\frac{\partial^{n} \Phi_{1}}{\partial c^{n}}\right)
= -p\,\frac{\partial^{n-1}\Phi_{1}}{\partial c^{n-1}}.
\]

For later reference it is convenient to write these relations in a compact way,
explicitly indicating the dependence on $n$:
\begin{equation}\label{eq:Phi1-D-rec}
	\begin{aligned}
		D\!\left(\frac{\partial^{n} \Phi_{1}}{\partial a^{n}}\right)
		&= n\,(x\,p + b)\,
		\frac{\partial^{n-1}\Phi_{1}}{\partial a^{n-1}},\\
		D\!\left(\frac{\partial^{n} \Phi_{1}}{\partial b^{n}}\right)
		&= n\,(x\,p + y\,q + a)\,
		\frac{\partial^{n-1}\Phi_{1}}{\partial b^{n-1}},\\
		D\!\left(\frac{\partial^{n} \Phi_{1}}{\partial c^{n}}\right)
		&= -n\,p\,
		\frac{\partial^{n-1}\Phi_{1}}{\partial c^{n-1}}.
	\end{aligned}
\end{equation}
Similarly, repeated differentiation of
\eqref{eq:first-a-M}–\eqref{eq:first-c-M} leads to
\begin{equation}\label{eq:Phi1-M-rec}
	\begin{aligned}
		M\!\left(\frac{\partial^{n} \Phi_{1}}{\partial a^{n}}\right)
		&= n\,\frac{\partial^{n-1}\Phi_{1}}{\partial a^{n-1}},\\
		M\!\left(\frac{\partial^{n} \Phi_{1}}{\partial b^{n}}\right)
		&= 0,\\
		M\!\left(\frac{\partial^{n} \Phi_{1}}{\partial c^{n}}\right)
		&= -n\,q\,
		\frac{\partial^{n-1}\Phi_{1}}{\partial c^{n-1}}.
	\end{aligned}
\end{equation}

Equations \eqref{eq:Phi1-D-rec} and \eqref{eq:Phi1-M-rec} constitute a system
of simple recurrence relations which can be used inductively to generate the
$n$th derivatives of $\Phi_{1}$ with respect to $a$, $b$ and $c$, once the
$(n-1)$st derivatives are known. When desired, the operators $(x\,p + b)$ and
$(x\,p + y\,q + a)$ acting on $\Phi_{1}$ or its derivatives can be replaced by
contiguous combinations of Humbert functions, using the differentiation
formulas in Section~\ref{sec:applications}.

\subsection{Other Humbert functions}

For the remaining Humbert confluent hypergeometric functions
$\Phi_{2},\Phi_{3},\Psi_{1},\Psi_{2},\Xi_{1},\Xi_{2}$ we follow exactly the same
strategy. Each of these functions satisfies a pair of second–order PDEs of
hypergeometric type. Denoting by $(P,Q)$, $(\widetilde{D},\widetilde{M})$,
$(D,M)$, etc., the corresponding pairs of differential operators, we again
differentiate the PDEs with respect to the parameters and obtain recursive
relations. Since the calculations are entirely analogous to those carried out
for $\Phi_{1}$, we only state the resulting formulas.

\subsubsection*{The function $\Phi_{2}$}

Let $P$ and $Q$ be the two PDE operators satisfied by $\Phi_{2}$. Then for
$n\ge 1$ we have
\begin{equation}\label{eq:Phi2-rec}
	\begin{aligned}
		P\!\left(\frac{\partial^{n}\Phi_{2}}{\partial a^{n}}\right)
		&= n\,\frac{\partial^{n-1}\Phi_{2}}{\partial a^{n-1}},\\
		P\!\left(\frac{\partial^{n}\Phi_{2}}{\partial b^{n}}\right)
		&= 0,\\
		P\!\left(\frac{\partial^{n}\Phi_{2}}{\partial c^{n}}\right)
		&= -n\,p\,
		\frac{\partial^{n-1}\Phi_{2}}{\partial c^{n-1}},
	\end{aligned}
\end{equation}
and
\begin{equation}\label{eq:Phi2-Q-rec}
	\begin{aligned}
		Q\!\left(\frac{\partial^{n}\Phi_{2}}{\partial a^{n}}\right)
		&= 0,\\
		Q\!\left(\frac{\partial^{n}\Phi_{2}}{\partial b^{n}}\right)
		&= n\,\frac{\partial^{n-1}\Phi_{2}}{\partial b^{n-1}},\\
		Q\!\left(\frac{\partial^{n}\Phi_{2}}{\partial c^{n}}\right)
		&= -n\,q\,
		\frac{\partial^{n-1}\Phi_{2}}{\partial c^{n-1}}.
	\end{aligned}
\end{equation}

\subsubsection*{The function $\Phi_{3}$}

For the function $\Phi_{3}(a;b;x,y)$, let $\widetilde{D}$ and $\widetilde{M}$
denote the corresponding PDE operators. Then
\begin{equation}\label{eq:Phi3-rec}
	\begin{aligned}
		\widetilde{D}\!\left(\frac{\partial^{n}\Phi_{3}}{\partial a^{n}}\right)
		&= n\,\frac{\partial^{n-1}\Phi_{3}}{\partial a^{n-1}},\\
		\widetilde{D}\!\left(\frac{\partial^{n}\Phi_{3}}{\partial b^{n}}\right)
		&= -n\,p\,
		\frac{\partial^{n-1}\Phi_{3}}{\partial b^{n-1}},
	\end{aligned}
\end{equation}
and
\begin{equation}\label{eq:Phi3-M-rec}
	\begin{aligned}
		\widetilde{M}\!\left(\frac{\partial^{n}\Phi_{3}}{\partial a^{n}}\right)
		&= 0,\\
		\widetilde{M}\!\left(\frac{\partial^{n}\Phi_{3}}{\partial b^{n}}\right)
		&= -n\,q\,
		\frac{\partial^{n-1}\Phi_{3}}{\partial b^{n-1}}.
	\end{aligned}
\end{equation}

\subsubsection*{The functions $\Psi_{1}$ and
	$\Psi_{2}$}

If we again use the notation $(D,M)$ for the pair of PDE operators of
$\Psi_{1}$ (the symbols are the same as for $\Phi_{1}$, but they act now on
$\Psi_{1}$), we obtain
\begin{equation}\label{eq:Psi1-rec}
	\begin{aligned}
		D\!\left(\frac{\partial^{n}\Psi_{1}}{\partial a^{n}}\right)
		&= n\,(x\,p + b)\,
		\frac{\partial^{n-1}\Psi_{1}}{\partial a^{n-1}},\\
		D\!\left(\frac{\partial^{n}\Psi_{1}}{\partial b^{n}}\right)
		&= n\,(x\,p + y\,q + a)\,
		\frac{\partial^{n-1}\Psi_{1}}{\partial b^{n-1}},\\
		D\!\left(\frac{\partial^{n}\Psi_{1}}{\partial c^{n}}\right)
		&= n\,p\,
		\frac{\partial^{n-1}\Psi_{1}}{\partial c^{n-1}},\\
		D\!\left(\frac{\partial^{n}\Psi_{1}}{\partial d^{n}}\right)
		&= 0,
	\end{aligned}
\end{equation}
and
\begin{equation}\label{eq:Psi1-M-rec}
	\begin{aligned}
		M\!\left(\frac{\partial^{n}\Psi_{1}}{\partial a^{n}}\right)
		&= n\,\frac{\partial^{n-1}\Psi_{1}}{\partial a^{n-1}},\\
		M\!\left(\frac{\partial^{n}\Psi_{1}}{\partial b^{n}}\right)
		&= 0,\\
		M\!\left(\frac{\partial^{n}\Psi_{1}}{\partial c^{n}}\right)
		&= 0,\\
		M\!\left(\frac{\partial^{n}\Psi_{1}}{\partial d^{n}}\right)
		&= -n\,q\,
		\frac{\partial^{n-1}\Psi_{1}}{\partial d^{n-1}}.
	\end{aligned}
\end{equation}

For $\Psi_{2}(a;b,c;x,y)$ we obtain analogous formulas:
\begin{equation}\label{eq:Psi2-D-rec}
	\begin{aligned}
		D\!\left(\frac{\partial^{n}\Psi_{2}}{\partial a^{n}}\right)
		&= n\,\frac{\partial^{n-1}\Psi_{2}}{\partial a^{n-1}},\\
		D\!\left(\frac{\partial^{n}\Psi_{2}}{\partial b^{n}}\right)
		&= -n\,p\,
		\frac{\partial^{n-1}\Psi_{2}}{\partial b^{n-1}},\\
		D\!\left(\frac{\partial^{n}\Psi_{2}}{\partial c^{n}}\right)
		&= 0,
	\end{aligned}
\end{equation}
and
\begin{equation}\label{eq:Psi2-M-rec}
	\begin{aligned}
		M\!\left(\frac{\partial^{n}\Psi_{2}}{\partial a^{n}}\right)
		&= n\,\frac{\partial^{n-1}\Psi_{2}}{\partial a^{n-1}},\\
		M\!\left(\frac{\partial^{n}\Psi_{2}}{\partial b^{n}}\right)
		&= 0,\\
		M\!\left(\frac{\partial^{n}\Psi_{2}}{\partial c^{n}}\right)
		&= -n\,q\,
		\frac{\partial^{n-1}\Psi_{2}}{\partial c^{n-1}}.
	\end{aligned}
\end{equation}

\subsubsection*{The functions $\Xi_{1}$ and
	$\Xi_{2}$}

Finally, for $\Xi_{1}(a,b,c;d;x,y)$, using again a suitable pair of PDE
operators $(D,M)$, we obtain
\begin{equation}\label{eq:Xi1-D-rec}
	\begin{aligned}
		D\!\left(\frac{\partial^{n}\Xi_{1}}{\partial a^{n}}\right)
		&= n\,(x\,p + c)\,
		\frac{\partial^{n-1}\Xi_{1}}{\partial a^{n-1}},\\
		D\!\left(\frac{\partial^{n}\Xi_{1}}{\partial b^{n}}\right)
		&= 0,\\
		D\!\left(\frac{\partial^{n}\Xi_{1}}{\partial c^{n}}\right)
		&= n\,(x\,p + a)\,
		\frac{\partial^{n-1}\Xi_{1}}{\partial c^{n-1}},\\
		D\!\left(\frac{\partial^{n}\Xi_{1}}{\partial d^{n}}\right)
		&= -n\,p\,
		\frac{\partial^{n-1}\Xi_{1}}{\partial d^{n-1}},
	\end{aligned}
\end{equation}
while
\begin{equation}\label{eq:Xi1-M-rec}
	\begin{aligned}
		M\!\left(\frac{\partial^{n}\Xi_{1}}{\partial a^{n}}\right)
		&= 0,\\
		M\!\left(\frac{\partial^{n}\Xi_{1}}{\partial b^{n}}\right)
		&= n\,\frac{\partial^{n-1}\Xi_{1}}{\partial b^{n-1}},\\
		M\!\left(\frac{\partial^{n}\Xi_{1}}{\partial c^{n}}\right)
		&= 0,\\
		M\!\left(\frac{\partial^{n}\Xi_{1}}{\partial d^{n}}\right)
		&= -n\,q\,
		\frac{\partial^{n-1}\Xi_{1}}{\partial d^{n-1}}.
	\end{aligned}
\end{equation}

For $\Xi_{2}(a,b;c;x,y)$ we similarly obtain
\begin{equation}\label{eq:Xi2-D-rec}
	\begin{aligned}
		D\!\left(\frac{\partial^{n}\Xi_{2}}{\partial a^{n}}\right)
		&= n\,(x\,p + b)\,
		\frac{\partial^{n-1}\Xi_{2}}{\partial a^{n-1}},\\
		D\!\left(\frac{\partial^{n}\Xi_{2}}{\partial b^{n}}\right)
		&= n\,(x\,p + a)\,
		\frac{\partial^{n-1}\Xi_{2}}{\partial b^{n-1}},\\
		D\!\left(\frac{\partial^{n}\Xi_{2}}{\partial c^{n}}\right)
		&= -n\,p\,
		\frac{\partial^{n-1}\Xi_{2}}{\partial c^{n-1}},
	\end{aligned}
\end{equation}
and
\begin{equation}\label{eq:Xi2-M-rec}
	\begin{aligned}
		M\!\left(\frac{\partial^{n}\Xi_{2}}{\partial a^{n}}\right)
		&= 0,\\
		M\!\left(\frac{\partial^{n}\Xi_{2}}{\partial b^{n}}\right)
		&= 0,\\
		M\!\left(\frac{\partial^{n}\Xi_{2}}{\partial c^{n}}\right)
		&= -n\,q\,
		\frac{\partial^{n-1}\Xi_{2}}{\partial c^{n-1}}.
	\end{aligned}
\end{equation}

\medskip

To summarise, all Humbert confluent hypergeometric functions admit simple and
parallel recurrence relations for derivatives of arbitrary order with respect
to their parameters. These relations are obtained by a uniform procedure based
on differentiating the underlying PDEs and are especially convenient when
combined with the contiguous–relation formulas of the next section, where the
operators $(x\,p + \lambda)$ and $(x\,p + y\,q + \lambda)$ are rewritten in
terms of Humbert functions with shifted parameters.

\section{Differentiation formulas for the Humbert confluent hypergeometric functions}
\label{sec:applications}

In this section we derive differentiation formulas with respect to the parameters
for the Humbert confluent hypergeometric functions of two variables. We also
obtain several reduction formulas with respect to the variables which express
higher–order derivatives in terms of contiguous Humbert functions. Throughout,
we continue to use the notation
\[
p := \frac{\partial}{\partial x},
\qquad
q := \frac{\partial}{\partial y},
\]
so that, for example, $x p$ stands for the differential operator
$x\,\partial/\partial x$.

The results obtained here may be viewed as explicit realizations of the operator
recurrence relations derived in Section~\ref{sec:nth-order}. For instance, in
the case of $\Phi_{1}$, the relations \eqref{eq:Phi1-D-rec}–\eqref{eq:Phi1-M-rec}
express the $n$th parameter derivatives in terms of the $(n-1)$st ones by means
of the differential operators $x p + b$, $x p + y q + a$ and $p$, acting on
$\Phi_{1}$ and its parameter derivatives. In the present section, we show that
these operators can be rewritten in a simple way as combinations of contiguous
Humbert functions with shifted parameters by acting on the double–series
definitions recalled in Section~\ref{sec:preliminaries}. In this manner, the
abstract operator recurrences of Section~\ref{sec:nth-order} are converted into
closed formulas for parameter derivatives and variable derivatives.

The proofs of the differentiation formulas follow a common pattern. One starts
from a series representation such as \eqref{eq:Phi1-def}–\eqref{eq:Xi2-def},
applies a simple shift identity for Pochhammer symbols (for example,
$(a+1)_{n} = (a)_{n}\bigl(1 + n/a\bigr)$), and then rearranges the resulting
series to recognize the defining series of a contiguous Humbert function. For
this reason we give a detailed proof only in a prototype case
(Theorem~\ref{thm:Phi1-contiguous}) and, for the remaining theorems, we merely
sketch or omit the proofs.

\begin{theorem}\label{thm:Phi1-contiguous}
	Let $\Phi_1 = \Phi_1(a,b;c;x,y)$ be the Humbert function defined in
	\eqref{eq:Phi1-def}. Then the following differentiation–contiguous relations
	hold:
	\begin{align}
		(x p + y q + a)\,\Phi_1(a,b;c;x,y)
		&= a\,\Phi_1(a+1,b;c;x,y), \label{eq:Phi1-a-shift}\\
		(x p + b)\,\Phi_1(a,b;c;x,y)
		&= b\,\Phi_1(a,b+1;c;x,y), \label{eq:Phi1-b-shift}\\
		(x p + y q + c-1)\,\Phi_1(a,b;c;x,y)
		&= (c-1)\,\Phi_1(a,b;c-1;x,y). \label{eq:Phi1-c-shift}
	\end{align}
	Consequently,
	\begin{equation}\label{eq:Phi1-contiguous-relation}
		(a - c + 1)\,\Phi_1(a,b;c;x,y)
		= a\,\Phi_1(a+1,b;c;x,y)
		- (c-1)\,\Phi_1(a,b;c-1;x,y).
	\end{equation}
\end{theorem}

\begin{proof}
	Using the identity
	\[
	(a+1)_{m+n} = (a)_{m+n}\left(1 + \frac{m+n}{a}\right)
	\]
	in the defining series of $\Phi_1(a+1,b;c;x,y)$ and comparing with the series
	for $\Phi_1(a,b;c;x,y)$, we obtain \eqref{eq:Phi1-a-shift} after a
	straightforward rearrangement. In the same way, using
	\[
	(b+1)_m = (b)_m\left(1 + \frac{m}{b}\right),
	\]
	we obtain \eqref{eq:Phi1-b-shift}. Finally, by means of
	\[
	\frac{1}{(c-1)_{m+n}}
	= \left(1 + \frac{m+n}{c-1}\right)\frac{1}{(c)_{m+n}},
	\]
	we arrive at \eqref{eq:Phi1-c-shift}. The relation
	\eqref{eq:Phi1-contiguous-relation} then follows by eliminating the operator
	$x p + y q$ from \eqref{eq:Phi1-a-shift} and \eqref{eq:Phi1-c-shift}. In view
	of \eqref{eq:Phi1-D-rec}, these formulas identify the operators occurring on
	the right–hand side of the recurrence relations of Section~\ref{sec:nth-order}
	with explicit shifts in the parameters.
\end{proof}

\begin{theorem}\label{thm:Phi1-mixed-param-deriv}
	The mixed parameter–derivative formulas
	\begin{align}
		(x p + y q + a)\,\frac{\partial}{\partial b}\Phi_1(a,b;c;x,y)
		&= a\,\frac{\partial}{\partial b}\Phi_1(a+1,b;c;x,y), \label{eq:Phi1-db}\\
		(x p + b)\,\frac{\partial}{\partial a}\Phi_1(a,b;c;x,y)
		&= b\,\frac{\partial}{\partial a}\Phi_1(a,b+1;c;x,y) \label{eq:Phi1-da}
	\end{align}
	hold. More generally, for every integer $n\ge 1$ we have
	\begin{align}
		(x p + y q + a)\,\frac{\partial^{n}}{\partial b^{n}}\Phi_1(a,b;c;x,y)
		&= a\,\frac{\partial^{n}}{\partial b^{n}}\Phi_1(a+1,b;c;x,y),
		\label{eq:Phi1-dnb}\\
		(x p + b)\,\frac{\partial^{n}}{\partial a^{n}}\Phi_1(a,b;c;x,y)
		&= b\,\frac{\partial^{n}}{\partial a^{n}}\Phi_1(a,b+1;c;x,y).
		\label{eq:Phi1-dna}
	\end{align}
\end{theorem}

\begin{proof}
	Differentiating \eqref{eq:Phi1-a-shift} and \eqref{eq:Phi1-b-shift} with
	respect to $b$ and $a$, respectively, and observing that $x p + y q + a$ and
	$x p + b$ do not depend on these parameters, we obtain
	\eqref{eq:Phi1-db}–\eqref{eq:Phi1-da}. Repeated differentiation yields
	\eqref{eq:Phi1-dnb} and \eqref{eq:Phi1-dna}, which are compatible with the
	recurrences \eqref{eq:Phi1-D-rec}–\eqref{eq:Phi1-M-rec}.
\end{proof}

\begin{theorem}\label{thm:Phi1-variable-deriv}
	For every integer $r\ge 1$, the derivatives of $\Phi_1$ with respect to the
	variables satisfy the reduction formulas
	\begin{align}
		\frac{\partial^{r}}{\partial x^{r}}
		\Phi_1(a,b;c;x,y)
		&= \frac{(a)_r (b)_r}{(c)_r}\,
		\Phi_1(a+r,b+r;c+r;x,y), \label{eq:Phi1-x-deriv}\\
		\frac{\partial^{r}}{\partial y^{r}}
		\Phi_1(a,b;c;x,y)
		&= \frac{(a)_r}{(c)_r}\,
		\Phi_1(a+r,b;c+r;x,y). \label{eq:Phi1-y-deriv}
	\end{align}
\end{theorem}

\begin{proof}
	Termwise differentiation of the defining double series of $\Phi_1$ with respect
	to $x$ gives
	\[
	\frac{\partial^{r}}{\partial x^{r}}\Phi_1(a,b;c;x,y)
	= \sum_{m,n\ge 0}
	\frac{(a)_{m+n}(b)_m}{(c)_{m+n}\,m!\,n!}\,
	(m)_r\,x^{m-r}y^n,
	\]
	where $(m)_r = m(m-1)\cdots(m-r+1)$ and the terms with $m<r$ vanish. Writing
	$(m)_r = \Gamma(m+1)/\Gamma(m+1-r)$ and shifting the summation index, we
	obtain \eqref{eq:Phi1-x-deriv}. The proof of \eqref{eq:Phi1-y-deriv} is
	similar and therefore omitted. When these relations are combined with
	\eqref{eq:Phi1-D-rec}–\eqref{eq:Phi1-M-rec}, they yield explicit expressions
	for all mixed derivatives with respect to parameters and variables.
\end{proof}


\begin{theorem}\label{thm:Phi2-contiguous}
	Let $\Phi_2 = \Phi_2(a,b;c;x,y)$ be the Humbert function defined in
	\eqref{eq:Phi2-def}. Then
	\begin{align}
		(x p + a)\,\Phi_2(a,b;c;x,y)
		&= a\,\Phi_2(a+1,b;c;x,y), \label{eq:Phi2-a-shift}\\
		(y q + b)\,\Phi_2(a,b;c;x,y)
		&= b\,\Phi_2(a,b+1;c;x,y), \label{eq:Phi2-b-shift}\\
		(x p + y q + c-1)\,\Phi_2(a,b;c;x,y)
		&= (c-1)\,\Phi_2(a,b;c-1;x,y). \label{eq:Phi2-c-shift}
	\end{align}
	In particular,
	\begin{equation}\label{eq:Phi2-contiguous-relation}
		(a + b - c + 1)\,\Phi_2(a,b;c;x,y)
		= a\,\Phi_2(a+1,b;c;x,y)
		+ b\,\Phi_2(a,b+1;c;x,y)
		- (c-1)\,\Phi_2(a,b;c-1;x,y).
	\end{equation}
\end{theorem}

\begin{theorem}\label{thm:Phi2-mixed-param-deriv}
	For every integer $n\ge 1$ we have
	\begin{align}
		(x p + a)\,\frac{\partial^{n}}{\partial b^{n}}\Phi_2(a,b;c;x,y)
		&= a\,\frac{\partial^{n}}{\partial b^{n}}\Phi_2(a+1,b;c;x,y),
		\label{eq:Phi2-dnb}\\
		(y q + b)\,\frac{\partial^{n}}{\partial a^{n}}\Phi_2(a,b;c;x,y)
		&= b\,\frac{\partial^{n}}{\partial a^{n}}\Phi_2(a,b+1;c;x,y).
		\label{eq:Phi2-dna}
	\end{align}
	These relations are the explicit counterparts of the recurrences
	\eqref{eq:Phi2-rec}–\eqref{eq:Phi2-Q-rec} for $\Phi_{2}$.
\end{theorem}

\begin{theorem}\label{thm:Phi2-variable-deriv}
	For every integer $r\ge 1$ the reduction formulas
	\begin{align}
		\frac{\partial^{r}}{\partial x^{r}}\Phi_2(a,b;c;x,y)
		&= \frac{(a)_r}{(c)_r}\,
		\Phi_2(a+r,b;c+r;x,y),
		\label{eq:Phi2-x-deriv}\\
		\frac{\partial^{r}}{\partial y^{r}}\Phi_2(a,b;c;x,y)
		&= \frac{(b)_r}{(c)_r}\,
		\Phi_2(a,b+r;c+r;x,y)
		\label{eq:Phi2-y-deriv}
	\end{align}
	hold.
\end{theorem}

\begin{theorem}\label{thm:Phi3-contiguous}
	Let $\Phi_3 = \Phi_3(a;b;x,y)$ be the Humbert function defined in
	\eqref{eq:Phi3-def}. Then
	\begin{align}
		(x p + a)\,\Phi_3(a;b;x,y)
		&= a\,\Phi_3(a+1;b;x,y),
		\label{eq:Phi3-a-shift}\\
		(x p + y q + b-1)\,\Phi_3(a;b;x,y)
		&= (b-1)\,\Phi_3(a;b-1;x,y).
		\label{eq:Phi3-b-shift}
	\end{align}
	Consequently,
	\begin{equation}\label{eq:Phi3-contiguous-relation}
		(a - b + 1)\,\Phi_3(a;b;x,y)
		= a\,\Phi_3(a+1;b;x,y)
		- (b-1)\,\Phi_3(a;b-1;x,y).
	\end{equation}
\end{theorem}

\begin{theorem}\label{thm:Phi3-mixed-param-deriv}
	For every integer $n\ge 1$,
	\begin{equation}\label{eq:Phi3-dnb}
		(x p + a)\,\frac{\partial^{n}}{\partial b^{n}}\Phi_3(a;b;x,y)
		= a\,\frac{\partial^{n}}{\partial b^{n}}\Phi_3(a+1;b;x,y).
	\end{equation}
\end{theorem}

\begin{theorem}\label{thm:Phi3-variable-deriv}
	For every integer $r\ge 1$ we have
	\begin{align}
		\frac{\partial^{r}}{\partial x^{r}}\Phi_3(a;b;x,y)
		&= \frac{(a)_r}{(b)_r}\,
		\Phi_3(a+r;b+r;x,y),
		\label{eq:Phi3-x-deriv}\\
		\frac{\partial^{r}}{\partial y^{r}}\Phi_3(a;b;x,y)
		&= \frac{1}{(b)_r}\,
		\Phi_3(a;b+r;x,y).
		\label{eq:Phi3-y-deriv}
	\end{align}
\end{theorem}


\begin{theorem}\label{thm:Psi1-contiguous}
	Let $\Psi_1 = \Psi_1(a,b;c,d;x,y)$ be the Humbert function of $\Psi$–type. Then
	\begin{align}
		(x p + y q + a)\,\Psi_1
		&= a\,\Psi_1(a+1,b;c,d;x,y), \label{eq:Psi1-a-shift}\\
		(x p + b)\,\Psi_1
		&= b\,\Psi_1(a,b+1;c,d;x,y), \label{eq:Psi1-b-shift}\\
		(x p + c-1)\,\Psi_1
		&= (c-1)\,\Psi_1(a,b;c-1,d;x,y), \label{eq:Psi1-c-shift}\\
		(y q + d-1)\,\Psi_1
		&= (d-1)\,\Psi_1(a,b;c,d-1;x,y). \label{eq:Psi1-d-shift}
	\end{align}
	These imply the contiguous relations
	\begin{align}
		(a - b - d + 1)\,\Psi_1
		&= a\,\Psi_1(a+1,b;c,d;x,y)
		- b\,\Psi_1(a,b+1;c,d;x,y) \notag\\
		&\quad - (d-1)\,\Psi_1(a,b;c,d-1;x,y), \label{eq:Psi1-contig-1}\\[0.5em]
		(a - c - d + 2)\,\Psi_1
		&= a\,\Psi_1(a+1,b;c,d;x,y)
		- (c-1)\,\Psi_1(a,b;c-1,d;x,y) \notag\\
		&\quad - (d-1)\,\Psi_1(a,b;c,d-1;x,y), \label{eq:Psi1-contig-2}\\[0.5em]
		(b - c + 1)\,\Psi_1
		&= b\,\Psi_1(a,b+1;c,d;x,y)
		- (c-1)\,\Psi_1(a,b;c-1,d;x,y). \label{eq:Psi1-contig-3}
	\end{align}
\end{theorem}

\begin{theorem}\label{thm:Psi1-mixed-param-deriv}
	For every integer $n\ge 1$ we have
	\begin{align}
		(x p + y q + a)\,\frac{\partial^{n}}{\partial b^{n}}\Psi_1
		&= a\,\frac{\partial^{n}}{\partial b^{n}}\Psi_1(a+1,b;c,d;x,y),
		\label{eq:Psi1-dnb}\\
		(x p + b)\,\frac{\partial^{n}}{\partial a^{n}}\Psi_1
		&= b\,\frac{\partial^{n}}{\partial a^{n}}\Psi_1(a,b+1;c,d;x,y).
		\label{eq:Psi1-dna}
	\end{align}
\end{theorem}

\begin{theorem}\label{thm:Psi1-variable-deriv}
	For every integer $r\ge 1$,
	\begin{align}
		\frac{\partial^{r}}{\partial x^{r}}\Psi_1(a,b;c,d;x,y)
		&= \frac{(a)_r (b)_r}{(c)_r}\,
		\Psi_1(a+r,b+r;c+r,d;x,y),
		\label{eq:Psi1-x-deriv}\\
		\frac{\partial^{r}}{\partial y^{r}}\Psi_1(a,b;c,d;x,y)
		&= \frac{(a)_r}{(d)_r}\,
		\Psi_1(a+r,b;c,d+r;x,y).
		\label{eq:Psi1-y-deriv}
	\end{align}
\end{theorem}


\begin{theorem}\label{thm:Psi2-contiguous}
	Let $\Psi_2 = \Psi_2(a;b,c;x,y)$ be the second Humbert function of $\Psi$–type.
	Then
	\begin{align}
		(x p + y q + a)\,\Psi_2
		&= a\,\Psi_2(a+1;b,c;x,y), \label{eq:Psi2-a-shift}\\
		(x p + b-1)\,\Psi_2
		&= (b-1)\,\Psi_2(a;b-1,c;x,y), \label{eq:Psi2-b-shift}\\
		(y q + c-1)\,\Psi_2
		&= (c-1)\,\Psi_2(a;b,c-1;x,y). \label{eq:Psi2-c-shift}
	\end{align}
	In particular,
	\begin{equation}\label{eq:Psi2-contig}
		(a - b - c + 2)\,\Psi_2
		= a\,\Psi_2(a+1;b,c;x,y)
		- (b-1)\,\Psi_2(a;b-1,c;x,y)
		- (c-1)\,\Psi_2(a;b,c-1;x,y).
	\end{equation}
\end{theorem}

\begin{theorem}\label{thm:Psi2-mixed-param-deriv}
	For every integer $n\ge 1$,
	\begin{equation}\label{eq:Psi2-dnb}
		(x p + y q + a)\,\frac{\partial^{n}}{\partial b^{n}}\Psi_2
		= a\,\frac{\partial^{n}}{\partial b^{n}}\Psi_2(a+1;b,c;x,y).
	\end{equation}
\end{theorem}

\begin{theorem}\label{thm:Psi2-variable-deriv}
	For every integer $r\ge 1$,
	\begin{align}
		\frac{\partial^{r}}{\partial x^{r}}\Psi_2(a;b,c;x,y)
		&= \frac{(a)_r}{(b)_r}\,
		\Psi_2(a+r;b+r,c;x,y),
		\label{eq:Psi2-x-deriv}\\
		\frac{\partial^{r}}{\partial y^{r}}\Psi_2(a;b,c;x,y)
		&= \frac{(a)_r}{(c)_r}\,
		\Psi_2(a+r;b,c+r;x,y).
		\label{eq:Psi2-y-deriv}
	\end{align}
\end{theorem}


\begin{theorem}\label{thm:Xi1-contiguous}
	Let $\Xi_1 = \Xi_1(a,b,c;d;x,y)$ be the Humbert function of $\Xi$–type. Then
	\begin{align}
		(x p + a)\,\Xi_1
		&= a\,\Xi_1(a+1,b,c;d;x,y), \label{eq:Xi1-a-shift}\\
		(y q + b)\,\Xi_1
		&= b\,\Xi_1(a,b+1,c;d;x,y), \label{eq:Xi1-b-shift}\\
		(x p + c)\,\Xi_1
		&= c\,\Xi_1(a,b,c+1;d;x,y), \label{eq:Xi1-c-shift}\\
		(x p + y q + d-1)\,\Xi_1
		&= (d-1)\,\Xi_1(a,b,c;d-1;x,y). \label{eq:Xi1-d-shift}
	\end{align}
	Consequently,
	\begin{align}
		(a - c)\,\Xi_1
		&= a\,\Xi_1(a+1,b,c;d;x,y)
		- c\,\Xi_1(a,b,c+1;d;x,y),
		\label{eq:Xi1-contig-1}\\[0.4em]
		(a + b - d + 1)\,\Xi_1
		&= a\,\Xi_1(a+1,b,c;d;x,y)
		+ b\,\Xi_1(a,b+1,c;d;x,y) \notag\\
		&\quad - (d-1)\,\Xi_1(a,b,c;d-1;x,y),
		\label{eq:Xi1-contig-2}\\[0.4em]
		(b + c - d + 1)\,\Xi_1
		&= b\,\Xi_1(a,b+1,c;d;x,y)
		+ c\,\Xi_1(a,b,c+1;d;x,y) \notag\\
		&\quad - (d-1)\,\Xi_1(a,b,c;d-1;x,y).
		\label{eq:Xi1-contig-3}
	\end{align}
\end{theorem}

\begin{theorem}\label{thm:Xi1-mixed-param-deriv}
	For every integer $n\ge 1$ the following differentiation formulas hold:
	\begin{align}
		(x p + a)\,\frac{\partial^{n}}{\partial b^{n}}\Xi_1
		&= a\,\frac{\partial^{n}}{\partial b^{n}}\Xi_1(a+1,b,c;d;x,y),
		\label{eq:Xi1-dnb}\\
		(x p + a)\,\frac{\partial^{n}}{\partial c^{n}}\Xi_1
		&= a\,\frac{\partial^{n}}{\partial c^{n}}\Xi_1(a+1,b,c;d;x,y),
		\label{eq:Xi1-dnc}\\
		(y q + b)\,\frac{\partial^{n}}{\partial a^{n}}\Xi_1
		&= b\,\frac{\partial^{n}}{\partial a^{n}}\Xi_1(a,b+1,c;d;x,y),
		\label{eq:Xi1-dna}\\
		(y q + b)\,\frac{\partial^{n}}{\partial c^{n}}\Xi_1
		&= b\,\frac{\partial^{n}}{\partial c^{n}}\Xi_1(a,b+1,c;d;x,y),
		\label{eq:Xi1-dnc2}\\
		(x p + c)\,\frac{\partial^{n}}{\partial a^{n}}\Xi_1
		&= c\,\frac{\partial^{n}}{\partial a^{n}}\Xi_1(a,b,c+1;d;x,y),
		\label{eq:Xi1-dna2}\\
		(x p + c)\,\frac{\partial^{n}}{\partial b^{n}}\Xi_1
		&= c\,\frac{\partial^{n}}{\partial b^{n}}\Xi_1(a,b,c+1;d;x,y).
		\label{eq:Xi1-dnb2}
	\end{align}
\end{theorem}

\begin{theorem}\label{thm:Xi1-variable-deriv}
	For every integer $r\ge 1$,
	\begin{align}
		\frac{\partial^{r}}{\partial x^{r}}\Xi_1(a,b,c;d;x,y)
		&= \frac{(a)_r (c)_r}{(d)_r}\,
		\Xi_1(a+r,b,c+r;d+r;x,y),
		\label{eq:Xi1-x-deriv}\\
		\frac{\partial^{r}}{\partial y^{r}}\Xi_1(a,b,c;d;x,y)
		&= \frac{(b)_r}{(d)_r}\,
		\Xi_1(a,b+r,c;d+r;x,y).
		\label{eq:Xi1-y-deriv}
	\end{align}
\end{theorem}


\begin{theorem}\label{thm:Xi2-contiguous}
	Let $\Xi_2 = \Xi_2(a,b;c;x,y)$ be the second Humbert function of $\Xi$–type.
	Then
	\begin{align}
		(x p + a)\,\Xi_2
		&= a\,\Xi_2(a+1,b;c;x,y), \label{eq:Xi2-a-shift}\\
		(x p + b)\,\Xi_2
		&= b\,\Xi_2(a,b+1;c;x,y), \label{eq:Xi2-b-shift}\\
		(x p + y q + c-1)\,\Xi_2
		&= (c-1)\,\Xi_2(a,b;c-1;x,y). \label{eq:Xi2-c-shift}
	\end{align}
	In particular,
	\begin{equation}\label{eq:Xi2-contig}
		(a - b)\,\Xi_2(a,b;c;x,y)
		= a\,\Xi_2(a+1,b;c;x,y)
		- b\,\Xi_2(a,b+1;c;x,y).
	\end{equation}
\end{theorem}

\begin{theorem}\label{thm:Xi2-mixed-param-deriv}
	For every integer $n\ge 1$,
	\begin{align}
		(x p + a)\,\frac{\partial^{n}}{\partial b^{n}}\Xi_2
		&= a\,\frac{\partial^{n}}{\partial b^{n}}\Xi_2(a+1,b;c;x,y),
		\label{eq:Xi2-dnb}\\
		(x p + b)\,\frac{\partial^{n}}{\partial a^{n}}\Xi_2
		&= b\,\frac{\partial^{n}}{\partial a^{n}}\Xi_2(a,b+1;c;x,y).
		\label{eq:Xi2-dna}
	\end{align}
\end{theorem}

\begin{theorem}\label{thm:Xi2-variable-deriv}
	For every integer $r\ge 1$,
	\begin{align}
		\frac{\partial^{r}}{\partial x^{r}}\Xi_2(a,b;c;x,y)
		&= \frac{(a)_r (b)_r}{(c)_r}\,
		\Xi_2(a+r,b+r;c+r;x,y),
		\label{eq:Xi2-x-deriv}\\
		\frac{\partial^{r}}{\partial y^{r}}\Xi_2(a,b;c;x,y)
		&= \frac{1}{(c)_r}\,
		\Xi_2(a,b;c+r;x,y).
		\label{eq:Xi2-y-deriv}
	\end{align}
\end{theorem}

\section{Applications}
\label{sec:apps}

The differentiation formulas obtained in the preceding sections can be used in a
variety of problems arising in physics, applied mathematics, engineering and
related areas. In particular, many applications require the evaluation of
Humbert confluent hypergeometric functions for parameter values that are close
to, but not exactly equal to, some reference set of numerator or denominator
parameters. In such situations, the explicit formulas for derivatives with
respect to the parameters provide a convenient tool for constructing local
parameter expansions and for analysing special parameter configurations.

\medskip

We first illustrate how the relations derived in Sections~2–4 simplify when the
parameters satisfy certain algebraic constraints.

\subsection*{Special parameter configurations}

Consider the Humbert function
\[
\Phi_1(a,b;c;x,y),
\]
whose defining series has been recalled in Section~\ref{sec:preliminaries}.
Along the diagonal $a=c$ in the parameter space, the function effectively
depends on a reduced number of parameters. Combining the first–order derivative
formulas with respect to $a$ and $c$ (see Section~2), one finds that
\begin{equation}\label{eq:Phi1-a-c-diagonal}
	\left.
	\left(
	\frac{\partial}{\partial a}
	+
	\frac{\partial}{\partial c}
	\right)
	\Phi_1(a,b;c;x,y)
	\right|_{a=c}
	= 0.
\end{equation}
Thus, when $a$ and $c$ are varied simultaneously along the diagonal $a=c$, the
value of $\Phi_1$ remains unchanged; in other words, $\Phi_1$ is locally
constant with respect to the combined variation of $a$ and $c$ subject to
$a=c$.

A completely analogous phenomenon occurs for the Humbert function
\[
\Psi_1(a,b;c,d;x,y).
\]
If we restrict to the diagonal $b=c$, then the first–order parameter derivative
formulas for $\Psi_1$ imply that
\begin{equation}\label{eq:Psi1-b-c-diagonal}
	\left.
	\left(
	\frac{\partial}{\partial b}
	+
	\frac{\partial}{\partial c}
	\right)
	\Psi_1(a,b;c,d;x,y)
	\right|_{b=c}
	= 0.
\end{equation}
Hence, in this case $\Psi_1$ is invariant under simultaneous variations of $b$
and $c$ along the line $b=c$ in the $(b,c)$–plane.

Relations of the type \eqref{eq:Phi1-a-c-diagonal} and
\eqref{eq:Psi1-b-c-diagonal} are typical for many problems in which the
physical or geometric model singles out particular combinations of numerator
and denominator parameters.

\subsection*{Taylor expansions with respect to parameters}

An important class of applications of the parameter–derivative formulas
obtained in Sections~2–4 consists of Taylor expansions with respect to the
parameters. Such expansions allow one, for instance, to approximate Humbert
confluent hypergeometric functions near a given set of parameter values, or to
study their sensitivity with respect to small perturbations of those
parameters.

Let us fix a reference value $a_0\in\mathbb{C}$ and expand $\Phi_1$ with
respect to the parameter $a$ around $a=a_0$. Using the existence of the
$n$th–order derivatives $\partial^n \Phi_1 / \partial a^n$ and the formulas
established in Section~3, we obtain the Taylor expansion
\begin{equation}\label{eq:Phi1-Taylor-a}
	\Phi_1(a,b;c;x,y)
	=
	\sum_{n=0}^{\infty}
	\frac{(a-a_0)^n}{n!}
	\left.
	\frac{\partial^n}{\partial a^n}
	\Phi_1(a,b;c;x,y)
	\right|_{a=a_0},
\end{equation}
whenever the series converges. In the same way, fixing a reference value
$c_0\in\mathbb{C}$ and expanding with respect to the denominator parameter $c$
gives
\begin{equation}\label{eq:Phi1-Taylor-c}
	\Phi_1(a,b;c;x,y)
	=
	\sum_{n=0}^{\infty}
	\frac{(c-c_0)^n}{n!}
	\left.
	\frac{\partial^n}{\partial c^n}
	\Phi_1(a,b;c;x,y)
	\right|_{c=c_0}.
\end{equation}

The explicit expressions for the derivatives with respect to $a$ and $c$,
obtained earlier in terms of Srivastava’s triple hypergeometric function
$F^{(3)}$, turn the formal Taylor series
\eqref{eq:Phi1-Taylor-a}–\eqref{eq:Phi1-Taylor-c} into computable expansions.
Similar Taylor expansions can be written with respect to any other numerator or
denominator parameter, for $\Phi_1$ as well as for the remaining Humbert
functions $\Phi_2$, $\Phi_3$, $\Psi_1$, $\Psi_2$, $\Xi_1$ and $\Xi_2$, by
making use of the corresponding $n$th–order parameter derivatives derived in
Section~3.

These parameter expansions, together with the reduction formulas for variable
derivatives obtained in Section~4, provide a flexible framework for analytical
and numerical investigations of Humbert confluent hypergeometric functions in a broad range of applications.

\section{Numerical illustrations}
\label{sec:numerical}

In this section we present some numerical examples and graphical
representations which illustrate the differentiation formulas obtained in the
preceding sections. For concreteness we focus on the Humbert confluent
hypergeometric function
\[
\Phi_{1}(a,b;c;x,y)
= \sum_{m,n=0}^{\infty}
\frac{(a)_{m+n} (b)_{m}}{(c)_{m+n}\, m!\, n!}\, x^{m} y^{n},
\]
together with its derivative with respect to the parameter \( a \). All
computations are carried out for the parameter choice
\[
a = \tfrac{3}{2}, \qquad b = \tfrac{3}{4}, \qquad c = \tfrac{5}{2},
\]
and for \((x,y)\) in the square \([0,0.8]\times[0,0.8]\), which lies well
inside the region of absolute convergence of the defining double series.

\subsection{Numerical evaluation}

For the numerical evaluation of \(\Phi_{1}\) we truncated the defining double
series to all terms with \(m+n\leq N_{\max}\), with \(N_{\max}=20\). In this
regime and for the parameter set specified above, the tail of the series is
very small, so that the resulting approximation is sufficiently accurate for
illustrative purposes. The derivative with respect to \(a\) was computed by
differentiating the series termwise and using
\[
\frac{\partial}{\partial a}(a)_{m+n}
= (a)_{m+n}\bigl[\Psi(a+m+n) - \Psi(a)\bigr],
\]
where \(\Psi\) denotes the digamma function. Inserting this into the double
series for \(\Phi_{1}\) yields a rapidly convergent double series for
\(\partial \Phi_{1}/\partial a\).

Table~\ref{tab:numerics-Phi1} displays representative numerical values of
\(\Phi_{1}(a,b;c;x,y)\) and \(\partial \Phi_{1}/\partial a(a,b;c;x,y)\) for
several points \((x,y)\). All values are rounded to six decimal places.

\begin{table}[ht]
	\centering
	\caption{Sample numerical values of
		\(\Phi_{1}(a,b;c;x,y)\) and
		\(\partial \Phi_{1}/\partial a(a,b;c;x,y)\)
		for \(a=1.5\), \(b=0.75\), \(c=2.5\).}
	\label{tab:numerics-Phi1}
	\vspace{0.5em}
	\begin{tabular}{c c c c}
		\hline
		\(x\) & \(y\) &
		\(\Phi_{1}(a,b;c;x,y)\) &
		\(\displaystyle\frac{\partial \Phi_{1}}{\partial a}(a,b;c;x,y)\) \\
		\hline
		0.1 & 0.1 & 1.113819 & 0.079604 \\
		0.1 & 0.5 & 1.430456 & 0.325267 \\
		0.3 & 0.3 & 1.409546 & 0.318193 \\
		0.5 & 0.5 & 1.851702 & 0.749493 \\
		0.7 & 0.1 & 1.720895 & 0.717261 \\
		0.7 & 0.5 & 2.258974 & 1.277738 \\
		\hline
	\end{tabular}
\end{table}

From Table~\ref{tab:numerics-Phi1} we see that, for the chosen parameter
values, both \(\Phi_{1}\) and \(\partial \Phi_{1}/\partial a\) increase
monotonically as either \(x\) or \(y\) increases. Moreover, the derivative
\(\partial \Phi_{1}/\partial a\) is strictly positive at all sample points,
which reflects the fact that \(\Phi_{1}\) is increasing in the parameter \(a\)
for this range of variables.

\subsection{Two–dimensional plots}

To visualize more clearly the dependence of \(\Phi_{1}\) and its parameter
derivative on the variable \(x\), we consider the one–parameter family
\[
x \longmapsto \Phi_{1}(a,b;c;x,y_{0}),
\qquad
x \longmapsto \frac{\partial \Phi_{1}}{\partial a}(a,b;c;x,y_{0}),
\]
with the same parameter values as above and with \(y_{0}=0.3\) fixed. Both
functions were evaluated on a uniform grid in the interval \(0\leq x\leq 0.8\).

\begin{figure}[htbp]
	\centering
\includegraphics[width=0.5\textwidth]{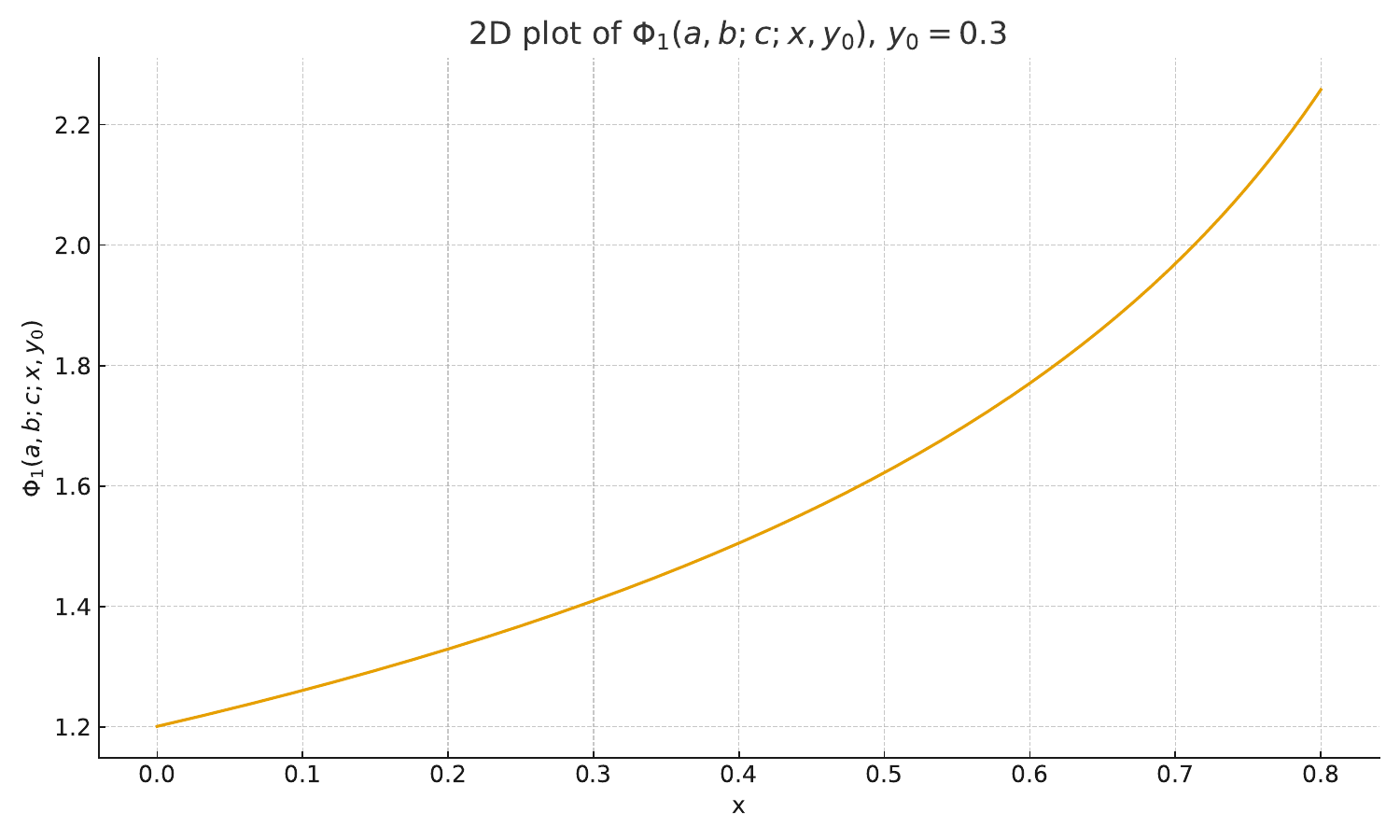}
	\caption{Plot of \(\Phi_{1}(a,b;c;x,y_{0})\) (solid line) and \(\partial \Phi_{1}/\partial a\) (dashed line) as functions
		of x for a=1.5\, b=0.75, c=2.5 and $y_{0}$=0.3.}
	\label{fig:Phi1-2D}
\end{figure}

Figure~\ref{fig:Phi1-2D} shows the resulting curves. The solid line corresponds
to \(\Phi_{1}(a,b;c;x,y_{0})\), while the dashed line corresponds to
\(\partial \Phi_{1}/\partial a\). As expected, both functions increase smoothly
with \(x\), and the derivative with respect to \(a\) grows more rapidly than
\(\Phi_{1}\) itself, indicating an enhanced sensitivity to changes in \(a\) as
\(x\) moves away from the origin.

\subsection{Three–dimensional surface plot}

We finally illustrate the joint dependence of \(\Phi_{1}\) on the variables
\((x,y)\) in the square \([0,0.8]\times[0,0.8]\). Using the same truncation
and parameter values as before, we computed \(\Phi_{1}(a,b;c;x,y)\) on a
uniform \(40\times40\) grid in this domain and constructed the corresponding
surface plot.

\begin{figure}[ht]
	\centering
	\includegraphics[width=0.7\textwidth]{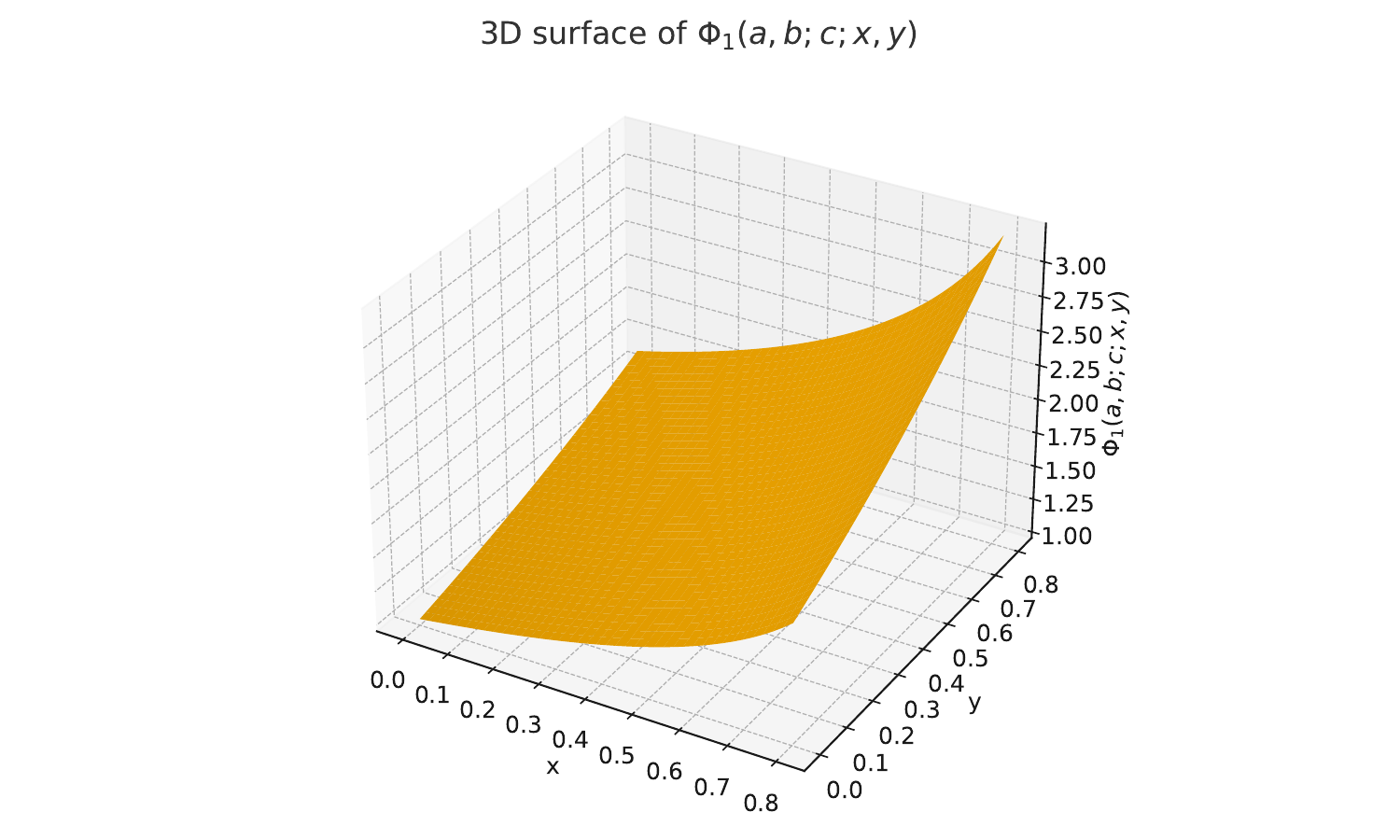}
	\caption{Three–dimensional surface plot of
		\(\Phi_{1}(a,b;c;x,y)\) for \(a=1.5\), \(b=0.75\), \(c=2.5\) and
		\((x,y)\in[0,0.8]\times[0,0.8]\).}
	\label{fig:Phi1-3D}
\end{figure}

The resulting graph is displayed in Figure~\ref{fig:Phi1-3D}. The surface is
smooth and strictly increasing in both variables, with a moderate curvature
near the origin and a steeper rise towards the corner \((x,y)=(0.8,0.8)\). This
behaviour is consistent with the positivity of the coefficients in the defining
double series of \(\Phi_{1}\) for the present choice of parameters.

These numerical examples provide a concrete illustration of the analytical
results obtained in the earlier sections, and they demonstrate that the
parameter–derivative formulas can be implemented efficiently in practical
computations.

\section{Concluding remarks}
\label{sec:conclusion}

In this text we have carried out a systematic study of derivatives with respect
to the parameters of the Humbert confluent hypergeometric functions of two
variables. More precisely, we considered all seven classical Humbert functions
\(\Phi_{1}\), \(\Phi_{2}\), \(\Phi_{3}\), \(\Psi_{1}\), \(\Psi_{2}\), \(\Xi_{1}\)
and \(\Xi_{2}\) and developed a unified framework for their differentiation
with respect to numerator and denominator parameters.

Starting from the double--series representations recalled in
Section~\ref{sec:preliminaries}, and using elementary properties of the Gamma,
digamma and polygamma functions, we first derived explicit formulas for the
first--order derivatives with respect to each parameter. These were then
recast in a compact and uniform manner in terms of Srivastava's triple
hypergeometric function \(F^{(3)}\), which plays a natural rôle as a basic
building block for multivariable parameter derivatives.

A second main ingredient of our approach is the use of the systems of linear
partial differential equations satisfied by the Humbert functions. By
differentiating these systems with respect to the parameters, we obtained
simple operator recurrences for parameter derivatives of arbitrary order, as
described in Section~\ref{sec:nth-order}. In Section~\ref{sec:applications}
these recurrences were combined with shift identities for Pochhammer symbols to
produce explicit differentiation and reduction formulas, expressing parameter
derivatives in terms of contiguous Humbert functions and higher--order
derivatives with respect to the variables.

To complement the theoretical developments, Section~\ref{sec:numerical}
presented a basic numerical illustration for the function \(\Phi_{1}\) and its
derivative with respect to a numerator parameter. Sample values and two-- and
three--dimensional plots were obtained directly from the double--series
representations, thereby demonstrating that the parameter--derivative formulas
can be implemented in a straightforward and numerically stable way.

The results obtained here provide an analytic toolkit for working with Humbert
confluent hypergeometric functions in contexts where parametric dependence is
essential, such as sensitivity analysis, perturbation methods and parameter
fitting in applied models. Several directions for further research remain open.
One natural extension is to consider generalized Humbert--type and related
multivariable hypergeometric functions and to derive analogous parameter
derivative formulas for them. Another direction is the development of dedicated
numerical algorithms which exploit the present formulas to compute Humbert
functions and their parameter derivatives efficiently over wider regions of the
parameter and variable space. We hope that the present text will serve as a
useful starting point for such investigations and for further applications in
mathematical physics and applied analysis.

	\section*{Affiliations}
\address{{\bf Ayman Shehata}: Department of Mathematics, Faculty of Science, Assiut University, Assiut 71516, Egypt.}\\
\email{aymanshehata@science.aun.edu.eg, drshehata2009@gmail.com, drshehata2006@yahoo.com}\\
{\bf ORCID ID:  0000-0001-9041-6752}\\

\address{{\bf Recep \c{S}ahin}: Department of Mathematics,  Faculty of Arts and Sciences, K\i r\i kkale University, 71450, K\i r\i kkale, Turkey.}\\
\email{recepsahin@kku.edu.tr}\\
{\bf ORCID ID: 0000-0001-5713-3830 }

\address{{\bf O\u{g}uz Ya\u{g}c\i}: Department of Mathematics,  Faculty of Arts and Sciences, K\i r\i kkale University, 71450, K\i r\i kkale, Turkey.}\\
\email{oguzyagci26@gmail.com, 1588151031@kku.edu.tr}\\
{\bf ORCID ID: 0000-0001-9902-8094 }

\address{{\bf Shimaa I. Moustafa}: Department of Mathematics, Faculty of Science, Assiut University, Assiut 71516, Egypt.}\\
\email{shimaa1362011@yahoo.com, shimaa$_{-}$m.@science.aun.edu.eg}\\
{\bf ORCID ID:  0000-0001-8589-9948}
	
\end{document}